 \documentclass[11pt]{amsart}
 \usepackage{latexsym}
 \usepackage{graphics, epsfig}
 \newtheorem{thm}{Theorem}[section]
 
 \newtheorem{lem}{Lemma}[section]
 \newtheorem{cor}{Corollary}[section]

 \newcommand{\pr}{\emph{Proof} : }

\begin{document}

\title{NON-COMPLEX SYMPLECTIC $4$-MANIFOLDS WITH $b_{2}^+ =1$}

\author{Jongil Park}

\address{Department of Mathematics, Konkuk University \\
         1 Hwayang-dong, Kwangjin-gu \\
         Seoul 143-701, Korea}

\email{jipark@konkuk.ac.kr}

\thanks{This work was supported by grant No. 1999-2-102-002-3 from
        the interdisciplinary research program of the KOSEF}

\date{March 1, 2001}

\subjclass{57R17, 57R57}

\keywords{Canonical class, rational, ruled, SW-invariant, symplectic}

\begin{abstract}
 In this paper we give a criterion whether a given minimal symplectic
 $4$-manifold with $b_{2}^+ =1$ having a torsion-free canonical class
 is rational or ruled. As a corollary, we confirm that most of
 homotopy elliptic surfaces
 $\{E(1)_{K} | K \ \mathrm{is \ a \ fibered \ knot \ in} \ S^3 \}$
 constructed  by R. Fintushel and \mbox{R. Stern} are minimal
 symplectic $4$-manifolds with $b_{2}^+ =1$ which do not admit a
 complex structure.
\end{abstract}

\maketitle

\section{Introduction}

\markboth{JONGIL PARK}
         {NON-COMPLEX SYMPLECTIC $4$-MANIFOLDS WITH $b_{2}^+ =1$}

 As an application of Seiberg-Witten theory to symplectic $4$-manifolds,
 C. Taubes proved that every minimal symplectic $4$-manifold with
 $b_{2}^{+}>1$ satisfies $K^2 \geq 0$ (\cite{t}).
 Here $K = -c_{1}(TX) \in H^{2}(X\!:\!{\bf Z})$ denotes the canonical
 class associated to a compatible almost complex structure on $X$.
 In the case when $b_{2}^+ =1$, the situation is quite different, i.e.
 there are examples satisfying $K^2 > 0,\, K^2 = 0$ and $K^2 < 0$
 respectively. Hence it is an interesting problem to classify  minimal
 symplectic $4$-manifolds with $b_{2}^+ =1$.
 Note that most known minimal symplectic $4$-manifolds with $b_{2}^+ =1$
 except Barlow surfaces and Dolgachev surfaces are rational or ruled,
 i.e. ${\bf CP}^2, S^2\times S^2$ or $S^2$-bundles over a Riemann
 surface. Regarding this, Liu proved that if $K^2 < 0$, then it is
 an irrational ruled surface (\cite{ms}).
 In this paper, when $K^2 \geq 0$, we give a criterion whether a given
 symplectic $4$-manifold with $b_{2}^+ =1$ having a torsion-free
 canonical class is rational or ruled. Explicitly, we prove

\begin{thm}
\label{main-1}
 Suppose $X$ is a closed minimal symplectic $4$-manifold with
 $b_{2}^+ =1$ such that its canonical class is a torsion-free
 class of non-negative square. Then $X$ is rational or
 ruled if and only if its Seiberg-Witten invariant
 \mbox{$SW_{X}^{\circ}$} vanishes.
\end{thm}

 R. Gompf constructed infinitely many minimal symplectic
 $4$-manifolds which do not admit a complex structure using a
 fiber sum technique (\cite{g}). All those manifolds have $b_{2}^+ >
 1$. In the case when $b_{2}^+ =1$, it has not been much known
 about minimal symplectic $4$-manifolds except complex surfaces.
 In Section~\ref{sec-3} we compute the Seiberg-Witten invariant
 $SW_{E(1)_{K}}^{\circ}$, obtained by a small generic perturbation
 of Seiberg-Witten equations, of a family of homotopy elliptic surfaces
 $\{E(1)_{K}\,| \, K \ \mathrm{is \ a \ fibered \ knot \ in} \ S^3 \}$
 constructed  by \mbox{R. Fintushel} and \mbox{R. Stern} in~\cite{fs} and,
 as a corollary, we conclude

\begin{cor}
\label{main-2}
 If $K$ is a fibered knot in $S^3$ whose Alexander polynomial is
 non-trivial and is different from that of any $(p,q)$-torus knot,
 then $E(1)_{K}$ is a simply connected minimal symplectic $4$-manifold
 which do not admit a complex structure.
\\
\end{cor}

\section{The Seiberg-Witten Invariant}
\label{sec-2}
 In this section we briefly review the Seiberg-Witten invariant of
 smooth $4$-manifolds (\cite{km}, \cite{m} for details).

 Let $X$ be a closed, oriented smooth $4$-manifold with $b_{2}^+ >0$
 and a fixed metric $g$, and let $L$ be a characteristic line bundle on $X$,
 i.e. $c_{1}(L)$ is an integral lift of  $w_{2}(X)$. This determines a
 $Spin^{c}$-structure on $X$ which induces a unique complex spinor
 bundle $W \cong W^{+} \oplus W^{-}$, where $W^{\pm}$ is
 the associated  $U(2)$-bundles on $X$ such that
 $\mathrm{det}(W^{\pm}) \cong L$.
 Note that the Levi-Civita connection on $TX$  together with a unitary
 connection $A$  on $L$ induces a connection  $\nabla_{A} :
 \Gamma(W^{+}) \rightarrow \Gamma(T^{\ast}X\otimes  W^{+})$.
 This connection,  followed  by  Clifford multiplication, induces a
 $Spin^{c}$-Dirac operator
 $D_{A} :  \Gamma(W^{+}) \rightarrow \Gamma(W^{-})$.
 Then, for each self-dual $2$-form
 $h \in \Omega^2_{+}(X\!:\!{\mathbf R})$,
 the following pair of equations  for a unitary connection $A$ on $L$
 and a section $\Psi$ of $\Gamma (W^{+})$ are called the
 {\em perturbed Seiberg-Witten equations}:
\begin{equation}
  (SW_{g,h})\left\{ \begin{array}{ll}
              D_{A}\Psi  & =\ \  0     \label{sw-eq}\\
              F_{A}^{+_{g}} & = \ \ i(\Psi \otimes \Psi^{\ast})_{0} + ih\, .
                  \end{array} \right.
\end{equation}
 Here $F_{A}^{+_{g}}$ is the self-dual part of the curvature of $A$ with
 respect to a metric $g$ on $X$ and $(\Psi  \otimes \Psi^{\ast})_{0}$ is the
 trace-free part of $(\Psi \otimes  \Psi^{\ast})$ which is interpreted as an
 endomorphism of $W^{+}$. The gauge group $\mathcal{G}:=Aut(L)\cong
 Map(X,S^{1})$ acts on the  space
 $\mathcal{A}_{X}(L) \times \Gamma (W^{+})$  by
 \[g\cdot (A,\Psi) =(g\circ A\circ g^{-1},g\cdot \Psi) \]
 Since the set of solutions is invariant under the action, it induces
 an orbit space, called the {\em Seiberg-Witten  moduli space}, denoted by
 $M_{X,g,h}(L)$, whose formal dimension is
 \[ \mathrm{dim}M_{X,g,h}(L)=\frac{1}{4}(c_{1}(L)^{2}-3\sigma(X)-2e(X))\]
 where $\sigma(X)$ is the signature of $X$ and $e(X)$ is the Euler
 characteristic of  $X$. Note that if $b_{2}^{+}(X) > 0$ and
 $M_{X,g,h}(L)\neq \phi$, then for a generic  metric $g$ and a generic
 self-dual $2$-form $h$ on $X$ the moduli space $M_{X,g,h}(L)$ contains
  no  reducible  solutions, so that it is a compact, smooth manifold of the
 given  dimension. Furthermore the moduli space $M_{X,g,h}(L)$ is orientable
 and  its orientation is determined by a choice of orientation on
 \mbox{det$(H^{0}(X;{\mathbf R})\oplus H^{1}(X;{\mathbf R})\oplus
 H^{2}_{+}(X;{\mathbf R}))$}.

\hspace*{-1.8em} {\bf Definition} The {\em Seiberg-Witten invariant}
 ({\em for brevity, SW-invariant}) for a smooth $4$-manifold $X$ with
 $b_{2}^+ >0$ is a function
 $SW_{X} : Spin^{c}(X) \rightarrow  {\mathbf Z}$  defined by
 \begin{equation}
 SW_{X}(L) = \left\{ \begin{array}{ll}
  \int_{M_{X,g,h}} \mu^{d_{L}}  &  \mathrm{if\ \ dim}M_{X,g,h}(L):=2d_{L}
    \geq 0\ \ \mathrm{and\ \ even} \\
  \ \ \ \ \ \ \ 0 &  \mathrm{otherwise} \, .
                       \end{array}
               \right.
\end{equation}
 Here $\mu \in H^{2}(M_{X,g,h}(L);{\bf Z})$ is the first Chern class of
 the $S^1$-bundle
 \[ \widetilde{M_{X,g,h}}(L) = \{\mathrm{solutions} (A,\Psi)\}/Aut^{0}(L)
    \longrightarrow M_{X,g,h}(L) \]
 where $Aut^{0}(L)$ consists of gauge transformations which are
 the identity on the fiber of $L$ over a fixed base point in $X$.
 Note that if $b_{2}^{+}(X) >1$, the Seiberg-Witten invariant, denoted by
 $SW_{X}= \sum SW_{X}(L)\cdot e^{L}$, is a diffeomorphic invariant,  i.e.
 $SW_{X}$ does not depend on the choice of a generic metric on $X$ and a
 generic perturbation of Seiberg-Witten equations.
 Furthermore, only finitely many $Spin^{c}$-structures on $X$ have
 a non-zero Seiberg-Witten  invariant.  We say that  the characteristic
 line bundle $L$, equivalently a  cohomology class
 $c_{1}(L) \in H^{2}(X;{\mathbf Z})$,  is a {\em SW-basic class} of $X$
 if $SW_{X}(L)\neq 0$.

 When $b_{2}^{+}(X) =1$, the SW-invariant $SW_{X}(L)$ defined  in $(2)$
 above  depends  not only on a metric $g$ but also on a self-dual
 $2$-form $h$. Because  of  this fact, there are several types of
 Seiberg-Witten invariants for a smooth  $4$-manifold with $b_{2}^+=1$
 depending on how to perturb Seiberg-Witten  equations. We introduce
 three  types of SW-invariants and investigate how they  are related.
 First we allow all metrics and self-dual $2$-forms to perturb
 Seiberg-Witten equations. Then the SW-invariant $SW_{X}(L)$ defined
 in $(2)$ above  has generically two values which are determined by the
 sign of $(2\pi  c_{1}(L) +[h])\cdot [\omega_{g}]$, where $\omega_{g}$
 is a  unique  $g$-self-dual harmonic $2$-form  of norm one lying
 in the positive component of $H^{2}_{+}(X;{\mathbf R})$.  We denote the
 SW-invariant for a generic  metric $g$ and a  self-dual $2$-form $h$
 satisfying  $(2\pi c_{1}(L) +[h])\cdot [\omega_{g}]>0$ by $SW_{X}^+(L)$
  and  denote the other one by $SW_{X}^-(L)$. Then the wall crossing
 formula  tells us the relation between  $SW_{X}^+(L)$ and $SW_{X}^-(L)$.

\begin{thm}[Wall crossing formula, \cite{ll}, \cite{m}]
\label{wall}
 Suppose that $X$ is a closed, oriented smooth $4$-manifold with
 $b_{1}=0$  and $b_{2}^+ =1$. Then for each characteristic line bundle
 $L$ on $X$ such  that the formal dimension of the moduli space
 $M_{X,g,h}(L)$ is non-negative  and even, say $2d_{L}$, we have
 \[ SW_{X}^+(L)  -  SW_{X}^-(L)  = -(-1)^{d_{L}} \, .\]
\end{thm}

\hspace{-1.8em} In general case ($b_{1}\neq 0$), Li and Liu  proved that
 the  difference is also expressed explicitly by some numerical
 invariants of $c_{1}(L)$ and $X$ (\cite{ll}).

  Second one may perturb the Seiberg-Witten equations by adding only a
 small generic self-dual $2$-form $h \in \Omega^2_{+}(X\!:\!{\mathbf R})$,
 so that one can define the SW-invariants as  in $(2)$ above. In this case
 we denote  the SW-invariant for a generic metric $g$ satisfying
 $(2\pi c_{1}(L))\cdot [\omega_{g}]>0$ by $SW_{X}^{\circ,+}(L)$
 and we denote the other one by  $SW_{X}^{\circ,-}(L)$.
  Note that $SW_{X}^{\circ,\pm}(L)=SW_{X}^{\pm}(L)$.
 But it sometimes happens that the sign of
 $(2\pi c_{1}(L))\cdot [\omega_{g}]$ is the same for all generic metrics,
 so that there exists only one SW-invariant obtained by a small generic
 perturbation of Seiberg-Witten equations.
 In such a case  we  define the SW-invariant  of $L$ on $X$ by
\begin{equation*}
 SW_{X}^{\circ}(L) = \left\{ \begin{array}{ll}
 SW_{X}^{\circ, +}(L)  &  \mathrm{if}\  2\pi c_{1}(L) \cdot [\omega_{g}] > 0 \\
 SW_{X}^{\circ, -}(L)  &   \mathrm{if}\  2\pi c_{1}(L) \cdot [\omega_{g}] < 0 \, .
                       \end{array}
               \right.
\end{equation*}
If $SW_{X}^{\circ}(L) \neq 0$,  we also call the corresponding $c_{1}(L)$
(or $L$) a {\em SW-basic class} which will be a diffeomorphic invariant of
 $X$.  Furthermore  we can extend many results obtained for smooth
 $4$-manifolds   with $b_{2}^+ >1$ to this case. For example,  we have

\begin{thm}
\label{basic-thm}
 Let $X$ be a  simply connected, closed, oriented smooth $4$-manifold with
 $b_{2}^+ =1$ and $b_{2}^- \leq 9$. Then  \\
 $($i$)$ There are only finitely many characteristic line bundles $L$ on $X$
         such that \mbox{$SW_{X}^{\circ}(L) \neq 0$}. \\
 $($ii$)$ If $X$ admits a metric of positive scalar curvature, then the
      SW-invariant of $X$ vanishes, that is, $SW_{X}^{\circ}(L) = 0$ for any
      characteristic line bundle $L$ on $X$. \\
 $($iii$)$ It also satisfies the generalized adjunction inequality. That is,
      if $\Sigma$ is a homologically non-trivial, smoothly embedded,
      oriented surface in $X$ with $[\Sigma]^{2}  \geq 0$,  then any
      characteristic line bundle $L$  on $X$ with \mbox{$SW_{X}^{\circ}(L)
      \neq 0$} satisfies the following inequality:
    $$2\cdot \mathrm{genus}(\Sigma) -2\ \geq \  |c_{1}(L) \cdot [\Sigma]|
      +  [\Sigma]^{2}\, . $$
\end{thm}

\pr\   Proofs of (i), (ii) and (iii) are exactly the same as the case
 $b_{2}^+ >1$ as long as the SW-invariant $SW_{X}^{\circ}$ is
 well defined,  i.e. it is independent of metrics on $X$.  Let $L$ be
 a  characteristic line  bundle on $X$ such that the formal dimension,
 $\frac{1}{4}(c_{1}(L)^2 - 3\sigma(X) -2e(X))$, of the moduli space is
 non-negative and even. The condition  $b_{2}^+ =1$ and $b_{2}^- \leq 9$
 imply  that  $c_{1}(L)^2 \geq 3\sigma(X) +2e(X) \geq 0$.
 Furthermore, since $X$ is simply  connected and $c_{1}(L)$ is
 characteristic,  $c_{1}(L) \neq 0$. Thus,  for any metric $g$ on $X$,
 the light  cone lemma (see Lemma~\ref{cone-lem}) implies
 $c_{1}(L)\cdot [\omega_{g}] \neq 0$,  so that the sign of
 $(2\pi c_{1}(L))\cdot [\omega_{g}]$ is the same for all  generic metrics.
 Hence the SW-invariant $SW_{X}^{\circ}(L)$ is well defined. $\ \ \ \Box$

 Finally we introduce one more type of Seiberg-Witten invariants for
 $b_{2}^+ =1$ - Given a fixed cohomology class
 $[x] \in H^{2}(X\!:\! {\mathbf Z})$ with
 $[x] \cdot [x] \geq 0$, one may divide a set of metrics and self-dual
 $2$-forms into two classes according to the sign of
 $proj_{+_{g}}(2\pi c_{1}(L) +[h])\cdot [x]$, where $proj_{+_{g}}$ is the
 projection of $\Omega^2(X\!:\!{\mathbf R})$ into  $g$-self-dual
 harmonic $2$-form space $H^2_{+_{g}}(X\!:\!{\mathbf R})$.
 In this case we denote the SW-invariant for a generic metric $g$
 and a generic self-dual \mbox{$2$-form} $h$ satisfying
 $proj_{+_{g}}(2\pi c_{1}(L)+[h])\cdot [x] >0$ by  $SW_{X}^{[x],+}(L)$
 and we denote the other one by $SW_{X}^{[x],-}(L)$.
  R. Fintushel and R. Stern used this type ($[x]=[T]$)  of SW-invariants
 for  $b_{2}^+ =1$ in \cite{fs}.  Note that if $[x] \cdot [\omega_{g}] > 0$,
 then  $SW_{X}^{[x],\pm}(L)= SW_{X}^{\pm}(L)$ and  if $[x] \cdot
 [\omega_{g}] < 0$, then $SW_{X}^{[x],\pm}(L)= SW_{X}^{\mp}(L)$.

 We close this section by mentioning that   C. Taubes'  theorem  on a
 symplectic $4$-manifold with $b_{2}^+ >1$  can be  easily extended
 to  $b_{2}^+ =1$  case.

\begin{thm}
\label{T-2}
 Suppose $X$ is a closed symplectic $4$-manifold with $b_{2}^+ =1$
 and a canonical class $K_{X}$. Then $SW_{X}^{-}(-K_{X}) = \pm 1$.
 \\
\end{thm}

\section{Minimal Symplectic $4$-Manifolds with $b_{2}^+ =1$}
\label{sec-3}
 In this section we investigate minimal symplectic $4$-manifolds with
 $b_{2}^+ =1$. First we compute the SW-invariant $SW_{X}^{\circ}$,
 obtained  by a small  generic perturbation of Seiberg-Witten  equations,
 of such symplectic $4$-manifolds.  And  then  we show that there exist
 infinitely many simply connected minimal symplectic $4$-manifolds which
 cannot admit a complex structure.
 Let us begin by stating  known results in the classification of symplectic
 $4$-manifolds. As mentioned in Introduction, most known minimal
 symplectic $4$-manifolds with $b_{2}^+ =1$ are complex surfaces
 such as rational or ruled surfaces, which are recognized by  Liu's
 and Ohta-Ono's results.

\begin{thm}[Liu, \cite{ms}]
\label{Liu-1}
 If $X$ is a closed  minimal symplectic $4$-manifold with $b_{2}^+ =1$
 such that  the canonical class $K_{X}$ has $K_{X}^2  < 0$,
 then $X$ is an irrational ruled  surface.
\end{thm}

 Furthermore,  when $K_{X}^2  \geq 0$, they  also proved the following fact
 by similar arguments used in the proof of Theorem~\ref{Liu-1} above.

\begin{thm}[Liu, Ohta-Ono, \cite{ms}]
\label{Liu-2}
  Let  $X$ be a closed  minimal symplectic $4$-manifold with $b_{2}^+ =1$
 and  a canonical class $K_{X}$. Then the followings are equivalent. \\
 $($i$)$  $X$ admits a metric of positive scalar curvature. \\
 $($ii$)$  $X$  admits a symplectic structure $\omega$ with
                   $K_{X}\cdot [\omega] < 0$. \\
 $($iii$)$  $X$ is either rational or ruled.
\end{thm}

 Hence, by combining two theorems above, we conclude that if $X$ is
 a minimal symplectic $4$-manifold  with $b_{2}^+ =1$
 which is neither  rational nor ruled,
 then it admits  a symplectic structure $\omega$ satisfying
 $K_{X}^2 \geq  0$ and  $K_{X}\cdot [\omega] \geq 0$.
 In this case  we compute the Seiberg-Witten invariant
 $SW_{X}^{\circ}(K_{X})$ of  the canonical class $K_{X}$ on $X$.

 Before stating  a main result, we introduce some terminology -
 Given a  symplectic $4$-manifold  $X$ with $b_{2}^+\! =\!1$ and
 a symplectic form  $\omega$,  the positive cone \mbox{$\{\alpha \in
 H^2(X\!:\!{\bf R})\,|\, \alpha^2 > 0 \}$} has two connected components
 in $H^2(X\!:\!{\bf R})$.  We write ${P}^+$ for the forward positive cone,
 which is the component containing $[\omega]$, and $P^-$ for the backward
 positive cone.  Then the  closure of ${P}^+$  in
$H^2(X\!:\!{\bf R})\setminus  \{0\}$ is given by
 \[ \overline{{P}^+ }= \{\alpha  \in  H^2(X\!:\!{\bf R}) \,|\,  \alpha^2
     \geq  0,\, \alpha \neq 0,\, \alpha \cdot [\omega] \geq 0 \} \, .\]
 Note that the next lemma, sometimes known as the light cone lemma,
 follows from Cauchy-Schwarz inequality.

\begin{lem}[\cite{ms}]
\label{cone-lem}
 If $\alpha, \beta \in \overline{P^+}$, then $\alpha \cdot \beta \geq 0$
 with equality if and only if $\alpha =\lambda \beta$ for some $\lambda >0$.
\end{lem}

 Now, in the case when $K^2 \geq 0$, we get a criterion whether a given
 symplectic $4$-manifold with $b_{2}^+ =1$ having a torsion-free
 canonical class is rational or ruled. Explicitly, we have

\begin{thm}
\label{main-thm1}
 Suppose $X$ is a closed minimal symplectic $4$-manifold with
 $b_{2}^+ =1$ such that its canonical class $K_{X}$ is a torsion-free
 class of non-negative square. Then $X$ is rational or ruled if
 and only if its Seiberg-Witten invariant \mbox{$SW_{X}^{\circ}$} vanishes.
\end{thm}

\pr\ As stated above, if $X$ is neither rational nor ruled,
 then it admits a symplectic structure $\omega$ satisfying
 $K_{X}^2 \geq  0$ and  $K_{X}\cdot [\omega] \geq 0$.  Furthermore,
 $K_{X}^2 = 3\sigma(X) + 2e(X) \geq 0$ implies $b_{2}^- \leq 9$
 and the light cone lemma implies  $K_{X}\cdot [\omega] > 0$.
 Note that, since  $b_{2}^- \leq 9$ and $(-K_{X})\cdot [\omega] < 0$,
 $SW_{X}^{\circ}$ is well defined and  $SW_{X}^{\circ}(-K_{X}) =
 SW_{X}^{\circ,-}(-K_{X}) = SW_{X}^{-}(-K_{X})$.
 Hence  Theorem~\ref{T-2} implies $SW_{X}^{\circ}(-K_{X})= \pm 1$,
 so that $SW_{X}^{\circ} \neq 0$.
 Conversely, if $X$ is rational or ruled with $K_{X}^2 \geq  0$,
 Theorem~\ref{Liu-2} implies that $X$ admits a metric of positive
 scalar curvature. Hence its SW-invariant $SW_{X}^{\circ}$
 vanishes by Theorem~\ref{basic-thm}. $\ \ \ \Box$\\

 Next, as an application of Theorem~\ref{main-thm1} above,
 we show that a family of homotopy elliptic surfaces
 $\{E(1)_{K}\,|\, K \ \mathrm{is \ a \ fibered \ knot \ in} \  S^3 \}$
 constructed by R. Fintushel and R. Stern in~\cite{fs} are minimal
 symplectic $4$-manifolds which are neither rational nor ruled.
 First we briefly  review their construction -
 Suppose $K$ is a  fibered knot  in  $S^{3}$ with a punctured surface
 $\Sigma_{g}^{\circ}$ of genus $g$  as fiber.
  Let $M_{K}$ be a $3$-manifold  obtained by performing $0$-framed
  surgery on  $K$, and let $m$ be a meridional circle to $K$.
  Then  the $3$-manifold  $M_{K}$ can be considered as a fiber bundle
 over circle with a closed Riemann surface $\Sigma_{g}$ as a fiber,
 and there is a  smoothly embedded torus $T_{m}:= m \times S^{1}$ of
 square $0$ in  $M_{K}\times S^{1}$. Thus $M_{K} \times S^{1}$ fibers
 over  $S^{1}\times S^{1}$ with $\Sigma_{g}$ as fiber and with
 $T_{m}= m \times  S^{1}$ as section.
 It is a Thurston's theorem that such a $4$-manifold $M_{K} \times S^{1}$
  has a symplectic structure with symplectic section $T_{m}$. Thus,
 if $X$ is a symplectic $4$-manifold with  a symplectically embedded
 torus  $T$ of square $0$, then the fiber sum  $4$-manifold
 \mbox{$X_{K}:= X\sharp_{T=T_{m}} (M_{K} \times S^{1})$},
 obtained by taking a  fiber sum along $T=T_{m}$, is symplectic.
 R. Fintushel and R. Stern  proved that  $X_{K}$ is homotopy equivalent
 to $X$  under a mild condition on  $X$  and computed the SW-invariant
 of $X_{K}$  (In the case when $b_{2}^+ =1$,
 they computed the relative SW-invariant of $X_{K}$).  For example,
 applying the construction above on an elliptic  surface $E(1)$,
 they get a   family of  homotopy  elliptic  surfaces
 $\{E(1)_{K}\,|\, K \ \mathrm{is \ a \ fibered \  knot \ in} \ S^3 \}$
 and computed  the relative SW-invariant  $SW^{\pm}_{E(1)_{K},T}$
 of $E(1)_{K}$:

\begin{thm}[\cite{fs}]
\label{fs-thm}
  For each fibered knot $K$ in $S^3$,  a homotopy elliptic surface
  $E(1)_{K}$  is a  simply connected symplectic $4$-manifold whose
 $[T]$-relative  SW-invariants are
 \begin{eqnarray*}
   \sum_{L\cdot [T] = 0} SW_{E(1)_{K}}^{[T],\pm}(L)\cdot e^L  & = &
   \sum_{L\cdot [T] = 0} SW_{E(1)}^{[T],\pm}(L)\cdot e^L \cdot
                         \Delta_{K}(e^{2[T]}) \\
    & = &  \sum^{\infty}_{n=0}(\mp 1)\cdot  e^{\mp(2n+1)[T]} \cdot
                         \Delta_{K}(e^{2[T]})
\end{eqnarray*}
 where $\Delta_{K}$ is the  Alexander polynomial of  $K$  and $T$ is
 a  symplectically embedded torus induced from a standard torus fiber
 lying in  $E(1)$.
\end{thm}

 Now we are ready to  compute the SW-invariant
 $SW_{E(1)_{K}}^{\circ}$ of $E(1)_{K}$, which is obtained  by
 a small generic perturbation of Seiberg-Witten equations.

\begin{thm}
\label{main-thm2}
 For each fibered knot $K$ in $S^3$, a homotopy elliptic surface $E(1)_{K}$
 has a SW-invariant denoted by $SW_{E(1)_{K}}^{\circ} =
 P_{E(1)_{K}}^{+} +P_{E(1)_{K}}^{-}$ such that  $P_{E(1)_{K}}^{\pm}$ is the
 partial sum consisting of only positive (negative) multiples of $[T]$
 in the exponent of the series
 $\mp \sum_{n=0}^{\infty} e^{\mp(2n+1)[T]} \cdot \Delta_{K}(e^{2[T]})$,
 where $\Delta_{K}$ is the Alexander polynomial of $K$ and $T$ is a
 symplectically embedded torus induced from a standard
 torus fiber lying in $E(1)$.
\end{thm}

\pr\   Since $E(1)_{K}$ is symplectic, one may choose a metric $g$
 on $E(1)_{K}$ so that the symplectic form $\omega$ is $g$-self-dual.
 Then a  symplectically  embedded torus $T$ of square $0$ in  $E(1)_{K}$,
 induced from a  standard torus fiber lying in $E(1)$, satisfies
 $[T]\cdot [\omega] > 0$.
 Thus   $ SW_{E(1)_{K}}^{\circ, \pm} = SW_{E(1)_{K}}^{\pm} =
 SW_{E(1)_{K}}^{[T],\pm}$.
  Since  $E(1)_{K}$ is  also simply connected and $b_{2}^- =9$,  each
 characteristic  line bundle $L$ on $E(1)_{K}$ has at most one SW-invariant,
 either $SW_{E(1)_{K}}^{+}(L)$ or $SW_{E(1)_{K}}^{-}(L)$. Hence
 $SW_{E(1)_{K}}^{\circ}(L)$ is well defined, i.e.
\begin{equation}
 SW_{E(1)_{K}}^{\circ}(L) = \left\{ \begin{array}{ll}
    SW_{E(1)_{K}}^{\circ, +}(L) = SW_{E(1)_{K}}^{[T], +}(L)  &
    \mathrm{if}\  2\pi c_{1}(L) \cdot [\omega] > 0 \\
   SW_{E(1)_{K}}^{\circ, -}(L) = SW_{E(1)_{K}}^{[T], -}(L)  &
    \mathrm{if}\  2\pi c_{1}(L) \cdot [\omega] < 0 \, .
                       \end{array}
               \right.
\end{equation}
 Furthermore  a generalized adjunction inequality
 (Theorem~\ref{basic-thm} - (iii)) implies that any  characteristic line
 bundle $L$ with   $SW_{E(1)_{K}}^{\circ}(L) \neq 0$ satisfies
 $L\cdot [T] =0$,  so that it is of the form $L=\lambda [T]$,
 for some integer $\lambda$,  due to the light cone lemma (see
 Lemma~\ref{cone-lem}).
 On the other hand,  Theorem~\ref{fs-thm} above implies that
 \begin{equation}
   \sum_{\lambda \in {\bf Z}} SW_{E(1)_{K}}^{[T], \pm}(\lambda [T])   =
   \mp \sum_{n=0}^{\infty} e^{\mp(2n+1)[T]} \cdot \Delta_{K}(e^{2[T]})
    \, .
\end{equation}
 Hence, combining Equations (3) and (4), we have
\begin{equation*}
   SW_{E(1)_{K}}^{\circ} :=  P_{E(1)_{K}}^{+} +P_{E(1)_{K}}^{-} \, ,
\end{equation*}
 where $P_{E(1)_{K}}^{\pm}$ is the  partial sum consisting  of
 only positive (negative) multiples of $[T]$ in  the exponent of the series
 $\mp \sum_{n=0}^{\infty} e^{\mp(2n+1)[T]} \cdot \Delta_{K}(e^{2[T]})$.
 $\ \ \ \Box$ \\

\begin{cor}
\label{main-cor2}
 For each fibered knot $K$ in $S^3$, a homotopy elliptic surface
 $E(1)_{K}$ is rational or ruled if and only if the Alexander
 polynomial $\Delta_{K}$ of $K$ is trivial.
\end{cor}

\pr\ Note that the \mbox{Alexander polynomial} $\Delta_{K}$ is
 non-trivial (i.e. $\Delta_{K} \neq 1$) if and only if the highest
 power of $[T]$ in the series
 $-\sum_{n=0}^{\infty} e^{-(2n+1)[T]} \cdot \Delta_{K}(e^{2[T]})$ coming
 from $(-e^{-[T]})\cdot$(the highest power of $\Delta_{K}(e^{2[T]})$) is
 non-negative, equivalently, the SW-invariant
 $SW_{E(1)_{K}}^{\circ}$ of $E(1)_{K}$ is not zero.
 Thus our claim follows from Theorem~\ref{main-thm1}. $\ \ \ \Box$ \\

\begin{cor}
\label{main-cor3}
 For each fibered knot $K$ with a non-trivial Alexander polynomial in $S^3$,
 a homotopy elliptic surface  $E(1)_{K}$ is a minimal symplectic $4$-manifold
 which is neither rational nor ruled.
\end{cor}

\pr\  Corollary~\ref{main-cor2} implies that $E(1)_{K}$ is not
 rational or ruled. Furthermore  since the set of SW-basic classes of
 $E(1)_{K}$ is $\{\lambda_{i}[T] \, | \, \lambda_{1}, \ldots , \lambda_{n} \,
  \mathrm{are \  some \  integers}\}$,
 we have  $(\lambda_{i}[T] - \lambda_{j}[T])^2 =0 \neq -4$.
 Hence  $E(1)_{K}$ is minimal.  $\ \ \ \Box$ \\

\hspace*{-1.8em} {\bf Examples.} Let $K$ be a $(p,q)-$torus knot
 with relatively prime integers $p, q \geq 2$ in $S^3$.
 Then $E(1)_{K}$ is diffeomorphic to a Dolgachev surface
 $E(1;p,q)$, which is obtained by $p-$ and $q-$ logarithmic transforms
 on $E(1)$. Hence the SW-invariant $SW_{E(1;p,q)}^{\circ}$ of
 Dolgachev surfaces is easily computed by Theorem~\ref{main-thm2}.
 Furthermore, it is an easy outcome of Corollary~\ref{main-cor3} that
 Dolgachev surfaces are neither rational nor ruled.

\begin{cor}
\label{main-cor4}
 If $K$ is a fibered knot in $S^3$ whose Alexander polynomial is
 non-trivial and is different from that of any $(p,q)$-torus knot,
 then $E(1)_{K}$ is a simply connected minimal symplectic $4$-manifold
 which do not admit a complex structure.
\end{cor}

\pr\ Note that simply connected complex surfaces satisfying
 $c_{1}^2 =0$ and $b_{2}^+ =1$ are classified as a blowing up of
 other complex surfaces and Dolgachev surfaces. Since the
 Seiberg-Witten invariant $SW_{E(1)_{K}}^{\circ}$ is different
 from that of any such complex surfaces, $E(1)_{K}$ is not
 diffeomorphic to any of them  so that it cannot be a complex
 surface. The rest of proof follows from Corollary~\ref{main-cor3}.
 $\ \ \ \Box$ \\

\hspace*{-1.8em} {\em Acknowledgement.} The author would like to
 thank Prof. Ron Fintushel for many helpful discussions
 including his explanation of Dolgachev surfaces. \\


\vspace*{2em}

\begin{thebibliography}{BDF}

\bibitem[FS]{fs} R. Fintushel and R. Stern,  Knots, links and 4-manifolds,
                 Invent. Math. {\bf 134} (1998), 363-400
\bibitem[G]{g}   R. Gompf, A new construction of symplectic manifolds,
                 Annals of Math. {\bf 142} no {\bf 3} (1995), 527-595
\bibitem[KM]{km} P.B. Kronheimer and T.S. Mrowka, The genus of embedded
                 surfaces in the projective plane, Math. Res. Letters
                 {\bf 1} (1994), 709-808
\bibitem[LL]{ll} T.J. Li and A. Liu, General wall crossing formula, Math.
                 Res. Letters {\bf 2} (1995), 797-810
\bibitem[M]{m}   J. Morgan, The Seiberg-Witten equations and applications
                 to the topology of smooth four-manifolds, Math. Notes
                 {\bf 44}, Princeton University Press, 1996
\bibitem[MS]{ms} D. McDuff and D. Salamon, A survey of symplectic
                 $4$-manifolds with $b_{2}^{+} =1$, Turkish Jour. Math.
                 {\bf 20} (1996), 47-61
\bibitem[T]{t}   C.H. Taubes, The Seiberg-Witten invariants and symplectic
                 forms, Math. Res. Letters {\bf 1} (1994), 809-822

\end{thebibliography}
\end{document}